\documentclass[12pt]{amsart}
\usepackage[dvips]{graphicx}
\newcommand{\F}{{\mathbb{F}}}

\newcommand{\Z}{{\mathbb{Z}}}

\newcommand{\Q}{{\mathbb{Q}}}
\newcommand{\R}{{\mathbb{R}}}

\newlength{\figboxwidth}
\def\defeq{\ensuremath{\stackrel{\mathrm{def}}{=}}}
\newcommand{\subgroup}{<}

\newcommand{\action}{\leadsto}
\newcommand{\arrow}{\rightarrow}

\newcommand{\compos}{\circ}

\newcommand{\bs}{\backslash}

\newcommand{\Aut}{\operatorname{Aut}}

\newcommand{\End}{\operatorname{End}}

\renewcommand{\mod}{\operatorname{mod}}

\newcommand{\PSL}{\operatorname{PSL}}

\newcommand{\SL}{\operatorname{SL}}
\newcommand{\GL}{\operatorname{GL}}

\newcommand{\Tr}{\operatorname{Tr}}

\newcommand{\Lk}{\operatorname{Lk}}

\newtheorem{theorem}{Theorem}[section]

\newtheorem{proposition}[theorem]{Proposition}
\newtheorem{lemma}[theorem]{Lemma}

\newtheorem{corollary}[theorem]{Corollary}
\newtheorem{definition}[theorem]{Definition}

\mathchardef\GG="321D
\input xy
\xyoption{all}
\numberwithin{equation}{section}
\newcommand{\Dir}{\operatorname{Dir}}
\newcommand{\Bat}{\operatorname{Bat}}

\newcommand{\Comm}{\operatorname{Comm}}
\newcommand{\CAT}{\operatorname{CAT}}
\newcommand{\hight}{\operatorname{ht}}
\newcommand{\proj}{\operatorname{pr}}
\newcommand{\Circ}{\operatorname{Circ}}
\newcommand{\PGL}{\operatorname{PGL}}
\newcommand{\Cayley}{\operatorname{Cay}}
\newcommand{\FT}{\mathfrak{T}}

\def\defeq{\ensuremath{\stackrel{\mathrm{def}}{=}}}
\title{Automata and square complexes.}
\author[Y. Glasner]{Yair Glasner}
\address{Department of Mathematics, University of Illinois at Chicago,
         851 South Morgan Street, Chicago, IL, 60607-7045, USA}
\email{yair@math.uic.edu}
\urladdr{http://www.math.uic.edu/~yair}
\author[S. Mozes]{Shahar Mozes}
\address{Institute of Mathematics, The Hebrew University                 \\
        Jerusalem 91904, Israel}
\email{mozes@math.huji.ac.il}

\date{}
\begin{document}
\bibliographystyle{amsalpha}
\begin{abstract}
    We introduce a new geometric tool for analyzing groups of finite automata. To each
    finite automaton we associate a square complex. The square complex is covered by
    a product of two trees iff the automaton is bi-reversible.
    Using this method we give examples of free groups and of Kazhdan groups
    which are generated by the different states of one
    finite (bi-reversible) automaton. We also reproduce the theorem of
    Macedo{\'n}ska, Nekrashevych, Sushchansky, on the connection
    between bi-reversible automata and the commensurator of a regular tree.
\end{abstract}
\maketitle

\section{Introduction}
\label{sed:introduction} If $X= \{x_1,x_2,\ldots,x_d\}$ is a finite set we denote by $X^{*}$
the free monoid generated by $X$ and by $F_X$ the free group generated by $X$.
\begin{definition}
A finite synchronous automaton is a quadruple of the form $A = \langle X,Q,\pi,\lambda
\rangle$, where:
\begin{itemize}
\item $X$ is a finite set, referred to as the alphabet.
\item $Q$ a finite set, referred to as the states.
\item $\pi:X \times Q \arrow Q$ is called the transition function.
\item $\lambda: X \times Q \arrow X$ is called the output function.
\end{itemize}
An initial automaton is just an automaton with a choice of a state $(A,q_0); \ q_0 \in Q$.
\end{definition}
\noindent The term finite means that the set of states is finite, and synchronous means that
the output function gives one alphabet letter rather then a general word in $X^{*}$. We will
not treat these more general automata, all our automata will be finite and synchronous.

An initial automaton $(A,q_0)$ defines a map $\FT(A,q_0):X^{*} \arrow X^{*}$ by the recursive
definition:
\begin{eqnarray*}
\FT(A,q_0)(\emptyset) & = & \emptyset          \\
\FT(A,q_0)(xw) & = & \lambda(x,q_0) \FT(A,\pi(x,q_0))w \qquad \forall x \in X, \ w \in X^{*}.
\end{eqnarray*}
This map from $X^{*}$ to itself does not preserve the monoid structure, it is merely a word
length preserving transformation of $X^{*}$ which is the same thing as a morphism of the
$d-$regular rooted tree. If this map is invertible for all $q \in Q$ we call the automaton $A$
invertible.
\begin{definition}
    Let $A$ be an automaton (resp. an invertible automaton)
    We denote by $S(A)$ (resp. $G(A)$) the semigroup (resp. group) of mappings of $X^{*}$
    generated by all different choices of initial automata $\{\FT(A,q)\}_{q\in Q}$.
\end{definition}

Recent years have seen some renewed interest in finite automata and the groups and semigroups
that they generate. Groups with interesting properties (e.g. infinite torsion groups, groups
with intermediate growth) were constructed from very simple automata \cite{}. On the other
hand well known groups (e.g. the lamplighter group) were realized as groups of type $G(A)$,
giving rise to new understanding of their properties \cite{GZ:Lamplighter,GLSZ:Atiyah_conj}.
We refer the readers to the beautiful paper \cite{GNS:Automata} and the references therein for
a survey on this subject.

Another theory that has seen many recent advances is the theory of finite square complexes and
in particular the theory of lattices acting on products of trees,
\cite{BM:Simple,BMZ:Representations,Wise:Phd,NR:CAT0}.

Our purpose in this paper is to introduce a new geometric tool to the theory of finite
automata. With each finite automaton we associate a square complex, conversely certain finite
square complexes give rise to finite automata.

Section \ref{sec:correspondance} is dedicated to describing the correspondence between square
complexes and automata and its basic properties. Some ideas and constructions, well known in
the theory of finite automata, become particularly appealing in this new geometric setting.
The idea of associating with every automaton its {\it{dual automaton}} by interchanging the
alphabet with the set of states is essentially due to Aleshin \cite{Aleshin:Free}. If the
automata is invertible we can also talk about the inverse and its dual and so on. In general,
up to $8$ different automata can be obtained by applying the inverse and dual operations to a
given one. From the geometric point of view the dual and inverse operations correspond to
certain geometric operations on square complexes, namely changing the vertical and horizontal
directions or inverting the vertical orientation respectively. An automaton is called
bi-reversible if all eight associated automata are well defined. Geometrically an automaton is
bi-reversible iff the universal covering of the associated square complex is a product of two
trees. We conclude section \ref{sec:correspondance} with a discussion of the beautiful paper
\cite{MNS:reversible_automata} in which Macedo{\'n}ska Nekrashevych and Sushchansky draw the
connection between bi-reversible automata and the commensurator of a regular (un-rooted) tree.
Similar ideas were developed independently in \cite{LMZ:superrigidity} where commensurator
elements are interpreted in terms of re-colorings of finite regular graphs.

An ongoing project is that of understanding the groups that can be generated by finite
automata. Somewhat surprisingly, for example, it is not so easy to generate a free group. In
\cite{Aleshin:Free} Aleshin gives an example of two initial automata that generate a free
group. Later examples of free groups generated by finite automata are given in
\cite{BS:GLn,Oliinyk:Free_products,Oliinyk:Free_groups}. In all these examples free groups (or
groups that contain free groups) are generated by different initial automata, so the
construction does not yield a group of the form $G(A)$. Sidki suggested
\cite{Sidki:One_rooted} one of the automata $A$ appearing in \cite{Aleshin:Free} as a
candidate for generating a free group $G(A)$.

Section \ref{sec:examples} is dedicated to the construction of new examples of automata. We
construct two families of bi-reversible automata, one giving rise to free groups and the other
to Kazhdan groups. For the free groups we construct automata $A$ such that $G(A)$ is free and
the different initial automata associated with $A$ form a symmetric set of free generators. We
do this by first constructing the associated square complexes; but our method is group
theoretic rather then combinatorial. We construct torsion free lattices in $\Aut(T_h) \times
\Aut(T_v)$ acting transitively on the vertices of the product of two regular trees $T_h \times
T_v$.

We wish to thank prof. Grigorchuk whose talk at Hebrew university inspired this paper.
We wish also to thank Laurent Bartholdi and Andrzej {\.Z}uk for many helpful comments and discussions.

\section{Automata and square complexes.}
\label{sec:correspondance} In this section we describe the square complexes associated with
finite automata and their basic properties. We begin by collecting basic facts and terminology
concerning automata and square complexes.
\subsection{Automata}
Finite synchronous automata were defined in the introduction. We define here invertible
automata and their inverses
\begin{definition}
\label{def:invertible} A synchronous automaton $A= (X,Q,\lambda,\pi)$ is called
{\it{invertible}} if $\lambda (\cdot,q): X \arrow X$ is a permutation for every $q \in Q$. If
$A$ is invertible we associate with it an inverse automaton $A^{-1} = \langle
X,\overline{Q},\overline{\lambda}, \overline{\pi} \rangle$ where $\overline{Q}$ is a set
isomorphic to $Q$ under the isomorphism $q \arrow \overline{q}$ and the output and transition
functions are given by the formulas $\overline{\lambda}(\cdot,\overline{q}) =
\lambda(\cdot,q)^{-1} \quad \forall q \in Q$ and $\overline{\pi}(x,\overline{q}) =
\overline{\pi(\overline{\lambda}(x,q),q)}$.
\end{definition}
\noindent {\it{Remark: }} It is easy to check that the maps $\FT(A,q):X^{*} \arrow X^{*}$ are
invertible for all $q \in Q$ iff the automaton $A$ is invertible in the sense of definition
\ref{def:invertible}. Moreover the inverse transformation is given by the inverse automaton
$(\FT(A,q))^{-1} = \FT(A^{-1},\overline{q}) \ \forall q \in Q$. Therefore our definition of
invertible automata coincides with the one given in the introduction.
\begin{definition}
    A morphism of synchronous automata
    $\Phi:\langle X,Q,\pi,\lambda \rangle \arrow \langle X',Q',\pi',\lambda' \rangle$
    is a pair of maps $\Phi_X:X \arrow X'$ and $\Phi_Q:Q \arrow Q'$ satisfying
    $\lambda'(\Phi_Q(q),\Phi_X(x)) = \Phi_X(\lambda (q,x))$ and
    $\pi'(\Phi_Q(q),\Phi_X(x)) = \Phi_Q(\pi(q,x))$.
\end{definition}
\begin{definition}
Let $A = \langle X,Q,\lambda,\pi \rangle$ be an automaton. The {\it{dual}} automaton is
obtained by exchanging between the output function and the transition function, the states of
the dual automaton will be the alphabet of the original and vice versa. Namely if $A = \langle
X,Q,\lambda,\pi \rangle$ then $A^{*} = \langle Q,X,\pi,\lambda \rangle.$
\end{definition}

\subsection{Square Complexes}
We recall some of the terminology pertaining to square complexes; all our notation is taken
from \cite{BM:Lattices}. A graph is a set of vertices $V$ and a set of edges $E$ with boundary
maps $o,t: E \arrow V$, called the {\it{origin and terminus maps}}, and a fixed point free
action $D_2 \cong \Z/2\Z \rightsquigarrow E$ called {\it{edge inversion}} which is compatible
with the boundary maps: $o\overline{e} = te, \ t\overline{e} = oe$. We denote by $\Circ_4$ the
``circle of length four graph'' with a set of vertices $\{1,2,3,4\}$ and set of edges
$\{[i,j]: i- j = \pm 1 \ (\mod 4)\}$.

A square complex $Y$ is given by a graph $Y^{(1)}$ (the one-skeleton), a set of squares $S$, a
boundary map associating with each square $s \in S$ a circle of length four: $\partial
s:\Circ_4 \arrow Y^{(1)}$ and a fixed point free action of the dihedral group of order $8$ on
the squares $D_4 \rightsquigarrow S$ which is compatible with the natural action of the same
group on the set of circles of length four. In other words: in the definition of a graph
(resp. of a square complex) we are keeping track of the two directed edges (resp. eight
directed squares) that correspond to each geometric edge (resp. face) and the group $D_2$
(resp. $D_4$) acts on these orientations.

We call a square complex VH if every vertex link is a bipartite graph (allowing for graphs
with multiple edges) and if there is a partition of the set of edges into vertical and
horizontal $E = V \amalg H$ which agrees with the bipartition of the graph on every link.
Note that our definition here is a little more general than that of \cite{BM:Lattices}: we
allow for VH-complexes even if the links are not complete bipartite graphs. A VH-complex will
be called a VH-T-complex if every link is a complete bipartite graph. An orientation on $E$ is
a choice of one directed representative for every geometric edge. A square complex will be
called {\it{directed}} if there is a choice of orientation on $E$ such that opposite edges of
each square are of the same orientation. Namely if $s \in S$ is a square then $\partial s
([1,2])$ and $\partial s ([4,3])$ have the same orientation as well as $\partial s([2,3])$ and
$\partial s ([1,4])$.

\subsection{The square complex associated with an automaton.}
\begin{definition}
Let $A = \langle X,Q,\lambda,\pi \rangle$ be a finite synchronous automaton. The {\it{square
complex associated with $A$}} denoted by $\Sigma (A)$ is a directed VH-square complex
$(V,E,S,\partial)$ such that.
\begin{itemize}
\item $V$ contains only one vertex $O$.
\item The $1-$skeleton is $E = X^{-} \amalg X^{+} \amalg Q^{-} \amalg Q^{+}$. Where $X^{\pm}$
      and $Q^{\pm}$ are two disjoint copies of $X$ and $Q$ respectively. The $\pm$ distinction
      describes the choice of orientation on $E$, namely
      $\overline{x^{+}} = x^{-} \ \forall x \in X$ and similarly for $Q$. The horizontal edges are
      $H = X^{\pm}$ and the vertical are $V=Q^{\pm}$.
\item The geometric squares are pairs $X \times Q$. With the boundary given by (see figure
\ref{fig:square}):
      $$ \partial (x,q) = x ^{+}\pi(q,x)^{+} \lambda(q,x)^{-} q^{-}.$$
      The eight directed squares corresponding to each geometric square are obtained by
      applying formally all possible symmetries of the square.
     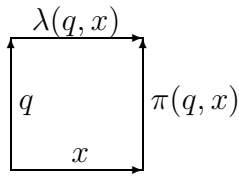
\begin{figure}
              \begin{center}
                \begin{picture}(80,80)
                  \put(10,10){\vector(1,0){50}}
                  \put(10,10){\vector(0,1){50}}
                  \put(10,60){\vector(1,0){50}}
                  \put(60,10){\vector(0,1){50}}
                  \put(33,13){\shortstack{$x$}}
                  \put(18,63){\shortstack{$\lambda(q,x)$}}
                  \put(13,33){\shortstack{$q$}}
                  \put(63,33){\shortstack{$\pi(q,x)$}}
                \end{picture}
                \caption{The square $(x,q)$} \label{fig:square}
              \end{center}
            \end{figure}
\end{itemize}
We denote by $\Pi(A) = \pi_1(\Sigma(A),O)$ and $\widetilde{\Sigma(A)}$ the fundamental group
and the universal cover of the VH-square complex $(V,E,S,\partial)$.
\end{definition}
\noindent {\it{Remark: }}
 A morphism of automata defines a morphism of labelled square complexes
$\Sigma(\Phi): \Sigma(A) \arrow \Sigma(A')$. The map $\Sigma(\Phi)$ in turn gives rise to maps
on the fundamental group and universal covering. We introduce the following short notation for
these maps $\tilde{\Phi} = \widetilde{\Sigma(\Phi)}:\widetilde{\Sigma(A)} \arrow
\widetilde{\Sigma(A)}$ and $\pi_1(\Phi)
\defeq \pi_1(\Sigma(\Phi)):\Pi(A) \arrow \Pi(A)$.

\subsection{Basic properties}
\begin{figure}
  \begin{center}
  \includegraphics{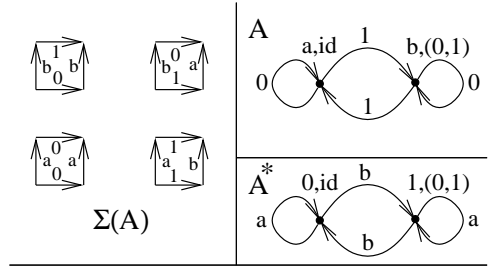}
  \end{center}
  \caption{A self-dual automaton defining the Lamplighter group.} \label{fig:lamplighter}
\end{figure}

Let $A$ be an automaton $\Sigma(A)$ the associated square complex. The group $D_4$ acts giving
rise to different choices of orientation and VH-structure on the same geometric complex. In
particular the complex associated with the dual automaton $\Sigma(A^{*})$ corresponds to the
square complex obtained by exchanging vertical and horizontal directions in $\Sigma(A)$. The
other six possibilities do not necessarily come form automata, in fact they come from the
following six automata iff the latter are well defined: $A^{-1}$, $(A^{-1})^{*}$,
$(A^{*})^{-1}$, $((A^{*})^{-1})^{*}$, $((A^{-1})^{*})^{-1}$, $(((A^{-1})^{*})^{-1})^{*} =
(((A^{*})^{-1})^{*})^{-1}$
\begin{definition} (See, \cite{MNS:reversible_automata})
An automaton is called reversible if its dual is invertible. It is called bi-reversible if all
eight automata above are well defined.
\end{definition}
\noindent {\it{Remarks: }}
\begin{itemize}
  \item It is enough to check that $A$ and $A^{-1}$ are reversible in order to
        deduce that $A$ is bi-reversible, i.e. it is enough that $A$,$A^{*}$ and
        $(A^{-1})^{*}$ all be invertible.
  \item It is possible for an automaton to be invertible and reversible but not bi-reversible.
        An example of such an automaton is depicted in figure \ref{fig:lamplighter}, it was shown
        in \cite{GZ:Lamplighter} (see also \cite{GNS:Automata}) that this automaton defines the
        lamplighter group.
\end{itemize}
\begin{proposition}
\label{prop:link_condition}
    Let $A$ be an automaton then A is bi-reversible iff $\Sigma(A)$ is
    VH-T-complex (i.e. iff $\widetilde{\Sigma(A)}$ is a product of trees).
\end{proposition}
\begin{proof}
    The link of the vertex $\Lk(O)$ is a bipartite graph where the vertices are naturally partitioned
    into two sets $X^{\pm}$ and $Q^{\pm}$. The orientation further divides each side
    into two disjoint sets. By the definition of $\Sigma(A)$,
    a pair of vertices of type $(x^{+},q^{+})$ are always connected by a unique edge in the
    link (i.e. there is a unique square incident with the two edges near their origin).
    The existence of such a unique edge for pairs of type $(x^{+},q^{-})$
    (resp. $(x^{-},q^{+})$,
    resp. $(x^{-},q^{-})$) is equivalent to the invertibility of $A$ (resp. $A^{*}$, resp.
    $(A^{-1})^{*})$. Thus $\Lk(O)$ is a complete bipartite
    graph if and only if the condition of the
    theorem holds. It is well known that the universal covering of a VH-complex is a product of
    trees iff every vertex link is a complete bipartite graph.
\end{proof}
\noindent
{\it{Remark: }}
    The condition of the theorem is also necessary and sufficient for $\widetilde{A}$ to be
    a $\CAT(0)$ space. Clearly a product of trees is $\CAT(0)$. Conversely, if
    the condition does not hold then, by counting considerations, there exists a pair of edges
    $(x,q)$ connected by more than one square. So there is a circle of length $\pi$ (two edges
    of length $\pi/2$ each) in the link
    $\Lk(O)$ and $\widetilde{\Sigma(A)}$ can not be $CAT(0)$.

Let $w \in X^{*}$, we can associate with it a path in the horizontal $1-$skeleton of the
complex $\Sigma(A)$. Such a path will follow always the positive orientation of the edges. A
general path in the horizontal $1-$skeleton is described by a word in $F_X$. In a similar
fashion we can describe paths in the vertical $1-skeleton$ by words in $Q^{*}$ or in $F_Q$. A
positive path in the $1-$skeleton of $\Sigma(A)$ is described by a word in $(X \amalg Q)^{*}$
and a general path by a word in $F_{X \amalg Q}$.

\begin{definition}
    An immersion of square complexes is a cellular map which induces an injection on the links.
\end{definition}
The fact that every pair of type $(x^{+},q^{+})$ is incident with a unique square is a direct
consequence of the definition of the complex $\Sigma(A)$. Using the notion of an immersion of
square complexes this can be extended to longer words in $X^{*}$ and $Q^{*}$. If one or more
of the automata $A,A^{*},(A^{-1})^{*}$ are invertible then we can talk about words in $F_{X}$
and $F_{Q}$ as well. The proof of the following theorem is a direct consequence of the
definitions.
\begin{proposition}
\label{prop:various_actions}
Let $A$ be an automaton.
\begin{itemize}
\item
   Let $\underline{x} \in X^{*}$ and $\underline{q} \in Q^{*}$, there is a unique immersion of
   a tessellated rectangle into $\Sigma(A)$
   whose lower edge describes the path $\underline{x}$ and whose left side describes the path
   $\underline{q}$.
\item
   If $A$ is invertible then the previous statement holds for every $\underline{x} \in X^{*}$ and
   $\underline{q} \in F_Q$.
\item
   If $A$ is reversible then one can take $\underline{x} \in F_{X}$ and
   $\underline{q} \in Q^{*}$.
\item
   If $A$ is bi-reversible then one can take
   $\underline{x} \in F_{X}$ and $\underline{q} \in F_{Q}$. Note that this is not possible if
   one only assumes that both $A,A^{*}$ are invertible.
\end{itemize}
\end{proposition}
The upper and right sides of the tessellated rectangle from the proposition describe elements
of these monoids (or groups, as the case may be: of the same type as the lower or left sides)
we will denote these by $\lambda(\underline{q},\underline{x})$ and
$\pi(\underline{q},\underline{x})$ respectively. This general notation agrees with our
previous notation in a few special cases:
 \begin{itemize}
 \item If $\underline{x} \in X, \underline{q} \in Q$ this agrees with the definition of ($\Sigma(A)$).
 \item If $q \in Q$ and $\underline{x} \in X^{*}$ then
       $\lambda(q,\underline{x}) = \FT(A,q)\underline{x}$.
 \item If $\underline{q} = q_1q_2 \ldots q_l$
       then $\lambda(\underline{q},\underline{x}) = \FT(A,q_1)\compos \FT(A,q_2) \compos \ldots
       \compos \FT(A,q_l)\underline{x}$. Thus the map:
       \begin{eqnarray*}
            Q^{*} \times X^{*} & \arrow & X^{*} \\
            (\underline{q} ,\underline{x}) & \arrow & \lambda(\underline{q},\underline{x})
       \end{eqnarray*}
       gives an action of $Q^{*}$ on $X^{*}$. This action is not necessarily effective, in fact,
       by definition of $S(A)$, it factors via the action of $S(A)$ on $X^{*}$.
 \item If $A$ is invertible (resp. reversible, resp. bi-reversible) we obtain more general actions
       of the form $F_Q \times X^{*} \arrow X^{*}$ (resp. $Q^{*} \times F_X \arrow F_X$ resp.
       $F_Q \times F_X \arrow F_X$). Where we remember that, by abuse of notation, the action on
       $F_X$ is not an action on the group but rather a word length preserving action or, what
       amounts to the same thing, an action on the Cayley graph $\Cayley(F_X,X)$ which fixes the
       base vertex $e$. In the invertible case the action $F_Q \times X^{*} \arrow X^{*}$ factors
       via an effective action of the group $G(A)$ on $X^{*}$.
\end{itemize}

\subsection{Non-finiteness properties}
Let $A_0$ be the automaton on a one letter alphabet and with one state. $\Sigma(A_0)$ is just
a torus obtained from gluing parallel sides of a square. We can identify $\Pi(A_0)
\rightsquigarrow \widetilde{\Sigma(A_0)}$ with the regular action of $\Z^{2}$ on $\R^{2}$.
Clearly there exists a unique morphism $\iota$ from any automaton $A$ to $A_0$. This gives
rise to {\it{height}} functions on the universal cover and fundamental groups of $A$. By abuse
of notation we call all of these functions $\hight = (\hight^{1},\hight^{2})$
\begin{equation}
\xymatrix{
       \Pi(A)  \ar[r]^{\hight}  \ar@{~>}[d]                   &
       \Z^{2}  \ar@{~>}[d]                              \\
       \widetilde{A} \ar[r]^{\hight}   &
       \R^{2}
}
\end{equation}
The following statement is a direct consequence of the existence of the height function:
\begin{lemma}
    Any path $w \in \Pi(A)$ homotopic to the identity satisfies $\hight(w) = (0,0)$.
    If $w$ is represented by a word in $F_{X \amalg Q}$ this means that in both horizontal and
    vertical directions the positively and negatively oriented edges should be balanced.
\end{lemma}
\begin{corollary}
    $\Pi(A)$ is infinite, in fact,
    every element $w \in \Pi(A)$ with $\hight(w) \ne (0,0)$ is of infinite order.
\end{corollary}

\subsection{Automata associated with square complexes}
As we have seen every square complex of the form $\Sigma(A)$ satisfies the minimal link
condition defined below.
\begin{definition}
A directed square complex is said to satisfy the {\it{minimal link condition}} if every pair
of positively oriented edges, one vertical and one horizontal edges starting at the same
vertex are incident to a unique square.
\end{definition}
Square complexes of the form $\Sigma = \Sigma(A)$ are all directed VH-complexes with one
vertex $O$ satisfying the minimal link condition. Conversely given such a square complex
$\Sigma$ we can reconstruct a finite automaton $A = A(\Sigma) = (X,Q,\pi,\lambda)$
\begin{itemize}
\item $X$ is the set of horizontal positively oriented edges.
\item $Q$ is the set of vertical positively oriented edges.
\item $\lambda(q,x)$ and $\pi(q,x)$ are defined to be the unique edges opposite to $x$ and $q$
      respectively, in the unique square incident to both $x$ and $q$ near their origin
      (see figure \ref{fig:square}).
\end{itemize}

If $\Sigma$ is a VH-T-complex then we can associate with it an automaton $A^{T}(\Sigma) =
(X,Q,\pi,\lambda)$ even without assuming that the complex is directed, by setting:
\begin{itemize}
\item $X$ is the set of directed horizontal edges.
\item $Q$ is the set of directed vertical edges.
\item $\lambda(q,x)$ and $\pi(q,x)$ are defined to be the unique edges opposite to $x$ and $q$
      in the unique square incident to $x$ and $q$ near their origin.
\end{itemize}
If we apply this construction to $\Sigma(A)$ where $A$ is a bi-reversible automaton then
we will end up with a new automaton whose sets of states and alphabet consist of the original
sets and their formal opposites.

Not every bi-reversible automaton is of the form $A^{T}(\Sigma)$. When it is possible,
however, there are big advantages for working with non directed square complexes. If $\Sigma$
is a directed square complex with universal covering $T_h \times T_v$ then the action of
$\pi_1(\Sigma,O)$ on each one of the factor trees can not be very transitive because of the
restriction that edge directions should be preserved. This is a serious drawback from the
point of view of the theory of groups acting on products of trees. For example in section
(\ref{sec:examples}) we construct bi-reversible automata from non-directed VH-T-square
complexes, using their strong transitivity properties in an essential way.

\begin{figure}
  \begin{center}
   \includegraphics{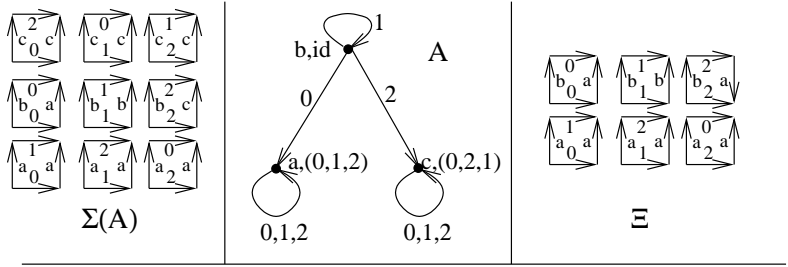}
  \end{center}
  \caption{An automata generating an infinite $3$ group.} \label{fig:pgroup}
\end{figure}

It is possible to formulate intermediate constructions when the link condition is defined in
two ``corners'' of the square complex (i.e. when $A(\Sigma)$ is invertible or when its dual is
invertible). Instead of trying to formulate a general statement we give an example. Consider
the automaton $A$ depicted in figure \ref{fig:pgroup}. This automaton is invertible but
neither $A$ nor $A^{-1}$ are reversible. It is proved in \cite{GNS:Automata} that $G(A)$ is
the Gupta-Sidki group which is an infinite $3$-group. Since $A$ has three states and acts on
an alphabet of three letters $\Sigma(A)$ is a complex of nine squares. If we notice however
that $\FT(A,a) = \FT(A,c)^{-1}$ we can draw a square complex $\Xi$ with six squares in such a
way that $A^{T}(\Xi)$ has four states $a,c=a^{-},b,b^{-}$ and the same alphabet, in such a way
that $\FT(A^{T}(\Xi),a) = \FT(A,a)$, $\FT(A^{T}(\Xi),a^{-}) = \FT(A,c)$, $\FT(A^{T}(\Xi),b) =
\FT(A,b)$, $\FT(A^{T}(\Xi),b^{-}) = \FT(A,b)^{-1}$.

\subsection{Bi-Reversible automata.}
From now on we will restrict our attention to VH-T-complexes $\Sigma$, and to bi-reversible
automata. As we have seen there are two cases here, the directed and non directed. We will
treat both. In the bi-reversible setting $\tilde{\Sigma} = T_h \times T_v$ is a product of
trees. The group $\pi_1(\Sigma,O)$ acts on $\tilde{A}$, preserving the VH-structure, and
therefore acts on each one of the factors.

Let $\tilde{O} = (O_h,O_v) \in T_h \times T_v$ be a base vertex. The horizontal and vertical
one skeletons $\Sigma^{(1),h},\Sigma^{(1),v}$ of the complex $\Sigma$ are just bouquets of
circles labelled by the set $X$ and $Q$ respectively. Since both $\Sigma^{(1),h}$ and $\Sigma$
are $\CAT(0)$ spaces then the embedding $\Sigma^{(1),h} \hookrightarrow \Sigma$ gives rise to
an embedding of the actions of the fundamental groups on the universal coverings see
\cite{BH:book}:
    $$\xymatrix{
       F_{X}  \ar[r]  \ar@{~>}[d]                   &
       \pi_1(\Sigma,O)  \ar@{~>}[d]                              \\
       T_h    \ar[r]   &
       T_h \times T_v}$$
And similarly in the vertical direction. We will henceforth identify the
trees $T_h \cong T_h \times \{O_v\}$ and $T_v \cong \{O_h\} \times T_v$  and the corresponding
groups with the setwise stabilizers of these trees $F_X \cong \pi_1(\Sigma,O)_{O_v}$ and
$F_Q \cong \pi_1(\Sigma,O)_{O_h}$.

The universal covering map to $\Sigma$ gives rise to an edge coloring of $T_h$ by the elements
of $X^{\pm}$. This gives a natural identification of $T_h$ with the Cayley graph of $F_X$
where $O_h$ represents the identity element of the group. The monoid $X^{*}$ spans an
$|X|$-regular rooted tree $T_X \subset T_h$.

\begin{figure}
  \begin{center}
  \includegraphics{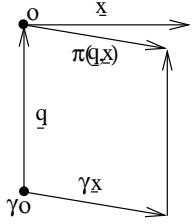}
  \end{center}
  \caption{The action of $\Pi(A)$ on $T_h$} \label{fig:action}
\end{figure}

Consider now the action of $F_Q \subset \Pi(A)$ on $T_h$. Let $q_1q_2 \ldots q_l =
\underline{q} \in F_Q$ and consider the vertical path $\gamma_{\underline{q}}$ in $\Sigma(A)$
represented by $\underline{q}^{-1} = q_l^{-1} \ldots q_2^{-1} q_1^{-1}$. This path represents
an element of $\pi_1(\Sigma,O)$ and it acts on $T_h \times T_v$ by a deck transformation.
$\gamma \tilde{O}$ is just the endpoint of the unique lifting $\tilde{\gamma}$ of the path
$\gamma$ starting at $\tilde{O}$. Let $\underline{x} = x_1x_2 \ldots x_m$ represent a
horizontal path in $\Sigma$ and $\tilde{\underline{x}}$ its unique lifting to $T_h \times T_v$
starting at $\tilde{O}$. The path $\gamma \tilde{\underline{x}}$ is the unique lifting of
$\underline{x}$ to a path starting at $\gamma \tilde{O}$. Finally projecting back to $T_h
\cong T_h \times \{O_v\}$ we find that $\proj_h(\gamma \tilde{\underline{x}})$ is the unique
lifting of the path $\lambda(\underline{q},\underline{x})$ to a path starting at $\tilde{O}$
(see figure \ref{fig:action}).
Summarizing everything we obtain the following
\begin{proposition}
\label{prop:action_description}
    Let $A = (X,q,\pi,\lambda)$ be a bi-reversible automaton, $\Sigma = \Sigma(A)$ the associated
    square complex and $T_h \times T_v$ its universal cover and $\tilde{O} = (O_h,O_v)$ a base
    vertex. The group $F_Q$ acts in two different ways on $\Cayley(F_X,X)$:
    \begin{itemize}
    \item Via the (left) action
          $(\underline{q},\underline{x}) \arrow \lambda(\underline{q},\underline{x})$
          defined after proposition \ref{prop:various_actions}.
    \item Via the identification $\pi_1(\Sigma,O)_{O_h} \cong F_Q$ and the (right)
          action of $\pi_1(\Sigma,O)$ on $T_h \times T_v$ which projects to an action on
          $T_h \cong \Cayley(F_X,X)$.
    \end{itemize}
    These two actions become isomorphic if we precompose the second action with the inverse map
    of the group in order to turn it into a left action.

    The rooted tree $X^{*}$ is invariant under these actions and any $g \in G(A)$ fixing
    pointwise $X^{*}$ fixes also the whole tree. The action of the group $F_Q$ on $T_h \cong
    \Cayley(F_X,X)$ is not necessarily effective but if factors through an effective action of
    the group $G(A)$.
\end{proposition}
\begin{proof}
    We have, in fact, proved everything with the exception of the last injectivity statement.
    Assume that $\underline{q} \in Q^{*}$ is such that
    $\lambda(\underline{q},\underline{x}) = \underline{x} \quad \forall \underline{x} \in X^{*}$.
    We first show that
    $\lambda(\underline{q},\underline{x}^{-1}) = \underline{x}^{-1} \quad \forall
    \underline{x} \in X^{*}$. Indeed consider the action of $\underline{q}$ on the periodic
    bi infinite
    word $w = \cdots \underline{x} \underline{x} \underline{x} \underline{x} \cdots$.
    In geometric terms we consider an immersion of an infinite strip whose
    lower side describes the word $w$ and some vertical segment describes the word $\underline{q}$.
    Since the possible number of vertical lines is finite and since the lower side and one
    vertical segment determine the immersion on the whole strip, the immersion is bound to be
    periodic. But by assumption the upper left side of the strip describes the word
    $\cdots \underline{x} \underline{x} \underline{x} \underline{x}$ so the upper right hand
    side must describe the word $\underline{x} \underline{x} \underline{x} \underline{x} \cdots$.
    The situation is completely symmetric of course: if we assume that $\underline{q}$ acts
    trivially on ``negative words'' we can conclude in a similar manner that it acts trivially
    on positive words as well.

    Clearly if the word $\underline{q}$ acts trivially on $X^{*}$ so does the word
    $\pi(\underline{q},\underline{x})$ for any $\underline{x} \in X^{*}$. Furthermore
    by the previous paragraph $\underline{q}$ will act trivially also on the set of negative
    words so that $\pi(\underline{q},\underline{x}^{-1})$ will also act trivially on all
    negative words and therefore $\pi(\underline{q},\underline{x}^{-1})$ will act trivially
    on $X^{*}$ as well for any $\underline{x} \in X^{*}$.

    Now decomposing a general word $w \in F_X$ into a sequence of negative and positive
    words we can prove by induction that $\lambda(\underline{q},w) = w$ and this completes the
    proof.
\end{proof}

The non directed situation is a little different: Let $\Sigma$ be a VH-T-complex with one
vertex, universal covering $T_h \times T_v$ and $A = A^{T}(\Sigma) = (X,Q,\pi,\lambda)$. Pick
an orientation on $X$. This gives a partition of $X$ into two sets $Y,\overline{Y}$. One can
think of $T_h$ as the Cayley graph of $F_Y$.  Let us define a homomorphism of monoids and of
trees.
 \begin{eqnarray*}
    \phi:X^{*} & \arrow & F_Y \\
         y     & \mapsto & y   \\
  \overline{y} & \mapsto & y ^{-1}
 \end{eqnarray*}
\begin{proposition}
\label{prop:action_description_T}
    The action $G(A) \action X^{*}$ factors, via $\phi$, to an action on $F_Y$.
    Furthermore the map $\Phi:G(A) \arrow \Aut(T_Y)$ is injective
    (i.e. any element acting trivially on $F_Y$ will act trivially also on $X^{*}$).
    The action of $G(A) \action F_Y$ is isomorphic to the action $F_Q \action T_h$ as above.
\end{proposition}
\begin{proof}
    Paths in the tree $\Cayley(F_y,Y)$ starting at the base vertex are in one to one
    correspondence with the set of paths without backtracking in the rooted tree
    corresponding to $X^{*}$ starting at the root. This gives a bijection between the
    automorphisms of $\Cayley(F_Y,Y)$ preserving the base point, and between the automorphisms
    of the rooted tree $\Cayley(X^{*},X)$ which preserve backtracking, i.e. automorphisms that
    take words of the type $wvv^{-1}$ to words of type $w'v'v^{'-1}$. The definition of $A$ implies
    $G(A)$ preserves ``backtracking''. Our discussion earlier
    proves that this correspondence of paths is equivariant under the $G(A)$ action.
\end{proof}
\begin{corollary}
In both the directed and the non directed cases  group $G(A)$ is isomorphic to the projection
of $F_Q$ to $\Aut(T_h)$.
\end{corollary}

\subsection{Bi-Reversible automata and the commensurator of a regular tree.}
Let $X$ be a finite set, $T = \Cayley(F_X,X)$ the right Cayley graph of $F_X$ with respect to
the free set of generators which is a $2|X|$-regular tree. The action of $F_X$ on $T$ gives
rise to an embedding of $F_X$ in $\Aut(T)$ as a uniform lattice. We define the commensurator
\begin{eqnarray*}
C &  = & \Comm_{\Aut(T)}(F_X) \\
  & = & \{g \in \Aut(T) | (gF_Xg^{-1} \cap F_X)
        {\text{ is of finite index in both }}F_X {\text{ and }} gF_Xg^{-1} \}.
\end{eqnarray*}
\begin{definition}
We call an automorphism of $T$ directed if it preserves the orientation on the edges which is
given by the choice of the set of generators. We denote the group of directed automorphisms of
$T$ by $\Aut^{\Dir}(T)$.
\end{definition}
Let $O \in T$ be a base vertex (e.g. corresponding to the identity element of $F_X$)
and $C_O$ the stabilizer of this point. Let $C^{\Dir}_O = C_O \cap \Aut^{\Dir}(T)$.

\noindent
{\it{Note: }} We have defined the commensurator for this specific lattice for convenience
but, in fact, a famous theorem of Leighton shows that up to conjugation in $\Aut(T)$
the commensurator is independent of the uniform lattice. This is not true anymore for the
directed commensurator.

We re-prove here a theorem from \cite{MNS:reversible_automata} using the new terminology.
\begin{theorem} (Macedo{\'n}ska, Nekrashevych and Sushchansky) \\
All bi-reversible automatic transformations over the alphabet $X$ form a group $\Bat(X)$. This
group is isomorphic to the group $C^{\Dir}_O$ defined above.
\end{theorem}
A non directed version of the same theorem will say:
\begin{theorem}
All automatic transformations of the form $(A^{T}(\Sigma),q)$ where $\Sigma$ is a VH-T-complex
with one vertex form a group. This group is naturally isomorphic to
a vertex stabilizer of the commensurator of a uniform lattice in $\Aut(T)$, where
$T$ is the $2|X|$-regular tree, i.e. to the group $C_O$ defined above.
\end{theorem}
\begin{proof}
    The proof of the two theorems is almost identical, we will treat the directed case
    and leave the other to the reader.
    It is enough to prove the isomorphism. First we define a natural homomorphism
    $\phi: \Bat(X) \arrow  C_O$.
    Let $T$ be the Cayley graph of $F_X$ which can also be identified as $T_h$: one of the factor
    trees of $\tilde{A}$ where $A$ is any finite automaton on the alphabet $X$.
    If $\FT(A,q) \in \Bat(X)$ is an initial bi-reversible automaton
    then, $\FT(A,q)$ is an automorphism of the rooted tree $T_X \subset T$ corresponding to
    the monoid $X^{*} \subset F_X$. However $\FT(A,q)$ can be extended to an automorphism
    $\widetilde{\FT(A,q)}$ of $T$. This can be thought of as
    the action of $\gamma(q^{-1}) \in \Pi(A)$, the loop corresponding to $q^{-}$ as an
    element of the fundamental group.  We define:
\begin{eqnarray}
\phi: \Bat(X) & \arrow & C_O   \nonumber \\
\phi      \FT(A,q) & \mapsto & \widetilde{\FT(A,q)}
\end{eqnarray}

    The map $\phi$ is well defined because a bi-reversible transformation $\FT(A,q)$ which is
    trivial on $T_X$ will give rise to a trivial $\widetilde{\FT(A,q)}$.

    To see that the image of $\phi$ is contained in the commensurator we identify
    $F_X \action \Cayley(F_X,X)$ with the action $\Pi(A)_{O_v} \action T_h \times \{O_v\}$.
    But $\gamma(q^{-1}) \Pi(A)_{O_v} \gamma(q) = \Pi(A)_{\gamma(q^{-1}) O_v}$ and these two
    groups are commensurable, namely the pointwise stabilizer of the geodesic
    $[O_2,\gamma(q^{-1})O_2]$ is a finite index intersection.

    Conversely assume that we are given an element $c \in C^{\Dir}_{O_1}$.
    By definition $(F_X \cap (c F_X c^{-1}))$ has a finite index in both $F_X$
    and $(c F_X c^{-1})$. By standard covering theory we obtain the corresponding
    covering maps of graphs:
$$\xymatrix{ & (F_X \cap (c F_X c^{-1})) \bs T \ar[dr]^{p_1}  \ar[dl]_{p_0}    & \\
               F_X \bs T  & &  c F_X c^{-1} \bs T             }
$$
    We take the product of the graph $Y = (F_X \cap (c F_X c^{-1})) \bs T$ with an interval
    $[0,1]$ and glue the graph ${0} \times Y$ onto  $F_X \bs T$ according to the covering
    map $p_o$ and the graph $\{1\} \times Y$ to $c F_X c^{-1} \bs T$ according to $p_1$.
    Finally we identify the two graphs $F_X \bs T$ and  $c F_X c^{-1} \bs T$ via the
    map $F_X e \arrow (c F_X c^{-1})c e$. It is easy to see, using the assumption on $c$,
    that the result is a directed VH-T-complex with one vertex describing the
    transformation $c$.
\end{proof}

\section{Examples}
\label{sec:examples}
\subsection{A free group}
Let $\Sigma$ be a VH-T-complex with one vertex and universal cover $T_h \times T_v$. Let $A =
A^{T}(\Sigma)$ be the associated automaton. By theorem \ref{prop:action_description_T} the
group $G(A)$ is the image in $\Aut(T_h)$ of the group $F_Q$---the subgroup of
$\pi_1(\Sigma,O)$ generated by all the vertical loops (i.e. loops labelled by elements of
$Q$). We recall that $F_Q$ is a free group with the set of states as free generators. The
image of $F_Q$ in $\Aut(T_h)$ will also be free if it acts faithfully on $T_h$. We are lead
examine the case where $\pi_1(\Sigma_O)$ acts faithfully on the factor trees of the universal
cover.
\begin{proposition}
Let $\Sigma$ be a VH-T-complex with one vertex and  $\pi_1(\Sigma,O) \action T_h \times T_v$
the action on the universal cover. If the action of $\pi_1(\Sigma,O)$ on the factor tree $T_h$
is effective then $G(A^{T}(\Sigma))$ is a free group on the set of states.
\end{proposition}
We proceed to develop the framework that enables us to ensure that this action is indeed
effective.
\begin{definition}
We say that a group action on a tree $G \action T$ is locally infinitely transitive,
if every vertex stabilizer $G_x$ acts transitively on each sphere
centered at $x$.
\end{definition}
\begin{proposition}
Let $\Sigma$ be a VH-T-complex, and $\pi_1(\Sigma,O) \leadsto T_h \times T_v$ its action on
the universal cover. If the action of $\pi_1(\Sigma,O)$ on $T_v$ is locally infinitely transitive
then its action on $T_h$ is effective.
\end{proposition}
\begin{proof}
    Assume the contrary and let
    $K \subgroup \pi_1(\Sigma,O)$ be the kernel of the action on $T_h$. The group
    $K$ is normal and it acts freely on $T_v$. The action of $\pi_1(\Sigma,O)$ on $T_v$
    factors now through an action of $K \bs \pi_1(\Sigma,O) \action K \bs T_h$.
    This new action should still be infinitely transitive (because dividing by
    a group acting freely does not change the stabilizers), but this is no longer
    possible for an action on a graph which is not a tree. This contradiction
    completes the proof.
\end{proof}
There are many ways of obtaining lattices acting locally infinitely transitively on the
projections (see \cite{BM:Lattices,BG:Affine}), but the easiest way is to appeal to the theory
of arithmetic lattices, in fact by the weak approximation theorem the projections of an
irreducible arithmetic lattice in $\PGL_2(\Q_p) \times \PGL_2(\Q_q)$ are dense in each of the
factor groups and in particular they act in a locally infinitely transitive fashion. Combining
these observations we have
\begin{proposition}
    Let $\Gamma < \PGL_2(\Q_p) \times \PGL_2(\Q_q)$ be a torsion free irreducible arithmetic
    lattice acting transitively on the vertices of the building $T_h \times T_v$,
    then $G^{T}(\Gamma \bs (T_h \times T_v))$ is a free group.
\end{proposition}
Examples of such lattices are described in detail in \cite{Vigneras:Quaternions,
moz:closures}. We will not go through the proof of the desired properties but will describe
the corresponding square complexes in detail for the convenience of the reader.

Assume that we are given a pair of primes $p,l$ both congruent to $1 (\mod (4))$. With such
a pair we associate a square complex:
\begin{itemize}
\item There is only one vertex $O$.
\item There are exactly $p+1$ integral quaternions $x=a+bi+cj+dk$ satisfying the
following:
    \begin{enumerate}
    \item $x = 1 (\mod 2)$ i.e. $a$ is odd and $b,c,d$ are all even.
    \item $Nx = x\overline{x} = a^2+b^2+c^2+d^2 = p$.
    \end{enumerate}
    Denote these by $x_1,x_2,\ldots,x_{p+1}$. Similarly there are $l+1$ quaternions
    $q_1,q_2,\ldots,q_{l+1}$ associated with the prime $l$.

    These quaternions will stand for the directed loops in the $1$ skeleton, where
    $x,\overline{x}$ are the two directed edges forming a geometric edge.
\item Given two quaternions $q_i,x_j$ there is a unique pair $q_{k(i,j)},x_{l(i,j)}$
      satisfying the equation:
    $$ x_j q_i = \pm q_{k(i,j)} x_{l(i,j)} \text{ (projective equality) }$$
      The set of these relations form the directed squares with boundary
      $x_j q_i \overline{q_{k(i,j)}} \overline{x_{l(i,j)}}$.
      Note that each non directed square appears here eight times in all possible orientations.
\end{itemize}

\subsection{A Kazhdan group}
\subsubsection{An extension of Mumford's example}
For the construction of an automaton $A$ such that $G(A)$ has property $T$ we will need to
construct an irreducible lattice $\Gamma < \PGL_3(\Q_p) \times \PGL_3(\Q_q)$ with the
following properties:
\begin{itemize}
 \item $\Gamma$ acts transitively on the vertices, of the building
       $\Delta_p \times \Delta_q$ of $\PGL_3(\Q_p) \times \PGL_3(\Q_q)$.
 \item $\Gamma$ acts freely on the geometric cubes of $\Delta_p^{(1)} \times \Delta_q^{(1)}$---
       the product of the one skeletons of the buildings.
\end{itemize}
Here $(p,q)$ are just a pair of prime numbers.

We start with an example due to Mumford (\cite{Mumford:Ample}) of a torsion free lattice
acting transitively on the building of $\PGL_3(\Q_2)$. Mumford's example is based on an
arithmetic lattice in $PGL_3(Q_2)$. For background on arithmetic lattices and on unitary
groups we refer to \cite{PR:Book} a good short summary on unitary groups in three variables
can be found in the first few pages of \cite{Rogawski:Unitary3}.

Let $\zeta = e^{2\pi i/7}$, $\lambda = \zeta + \zeta^{2} + \zeta^{4} = \tfrac{-1 +
\sqrt{-7}}{2}$ and $\overline{\lambda} = \zeta^{3} + \zeta^{5} + \zeta^{6} = \tfrac{-1 -
\sqrt{-7}}{2}$. Let $K = \Q(\lambda)$ and $\mathcal{O} = \Z[\lambda]$ its ring of integers.
The extension $[K:\Q]$ is of degree $2$ with a Galois involution --- the complex conjugation
$z \arrow \overline{z}$. The extension $[\Q(\zeta):K]$ is Galois of degree $3$, with the
Galois group generated by $\sigma$, $\sigma(\zeta) = \zeta^{2}$. We will, most of the time,
think of $\Q(\zeta)$ merely as a $K$-vector space and as such we will denote it by $V =
\Q(\zeta)$. $V$ contains the lattice $L = \Z[\zeta]$ with a basis $1,\zeta,\zeta^{2}$ over
$\mathcal{O}$. We define a Hermitian form on $V$:
$$
  h(x,y) = \Tr_{\Q(\zeta)/\Q(\lambda)}(x \overline{y})
         = x \overline{y} + \sigma(x \overline{y}) + \sigma^{2}(x \overline{y})
$$
Taking $1,\zeta,\zeta^{2}$ as a basis we find that $h(x,y) = xH(^{t}\overline{y})$ where $H$
is the matrix
$$
H = \begin{pmatrix}3 & \overline{\lambda} & \overline{\lambda}  \\
           \lambda & 3 & \overline{\lambda}         \\
           \lambda & \lambda & 3            \end{pmatrix}.
$$
Let $A = \End_{K}(V) = M_3(K)$ be the split central simple algebra over $K$. The form $h$
defines an involution of the second kind on $A$ given by $g \arrow \alpha(g) =
H(^{t}\overline{g})H^{-1}$. Let $G=U(h)$ be the algebraic group of all linear transformations
preserving the Hermitian form $h$. Note that, even though $A$ is a $K$-algebra, the algebraic
group $G$ is only defined over $\Q$ because the definition involves the Galois involution. The
$\Q$-rational points of $G$ are given by:
\begin{eqnarray*}
G(\Q) & = & \{g \in GL(V,K):h(xg,yg) = h(x,y) \ \forall x,y \in V \}.     \\
      & = & \{g \in A | g \alpha(g) = 1\}.
\end{eqnarray*}

Let $p$ be a place of $\Q$ (possibly $p = \infty$). We denote the $\Q_p$-rational points of
$G$ by
$$G_p = \{g \in M_3(K) \otimes_{\Q} \Q_p | g H (^t\overline{g}) H^{-1} = 1\}.$$
If $p = \mathfrak{P \overline{P}}$ splits in $K$ then $K_{\mathfrak{P}} \cong \Q_p$ and $K
\otimes_{\Q} \Q_p = K_{\mathfrak{P}} \oplus K_{\mathfrak{\overline{P}}} = \Q_p \oplus \Q_p$.
The natural extension of the Galois involution, acting trivially on $\Q_p$, switches the two
summands of $K \otimes_{\Q} \Q_p$. Consequently:
\begin{eqnarray*}
A_p & \defeq & A \otimes_{\Q} \Q_p \\
    & = & A_{\mathfrak{P}} \oplus A_{\mathfrak{\overline{P}}}.
\end{eqnarray*}
Up to isomorphism $\alpha$ induces the anti-involution $\alpha(g,h) = (^t h ,^t g)$ on
$A_p =  A_{\mathfrak{P}} \oplus A_{\mathfrak{\overline{P}}}$.
The projection into
each of the summands gives an isomorphism of $G_p$ with $\GL_3(K_{\mathfrak{P}}) \cong
\GL_3(\Q_p)$. Since the discriminant of $K$ is $d_K = -7$, the primes that split at $K$ are
exactly $2$ and the odd primes $p \ne 7$ such that $(\tfrac{-7}{p}) = 1$ (i.e.  all primes
congruent to $1,2$ or $4$ ($\mod 7$)).

Since $G_{\infty}$ is a compact group, by reduction theory if $p_1,p_2,\ldots,p_l$ are
rational primes all congruent to $1,2,4$ ($\mod 7$) then the group
\begin{eqnarray*}
\Gamma_1(p_1,\ldots,p_l) & = & G(\Z[1/p_1,1/p_2,\ldots,1/p_l])      \\
         & = & Q(\lambda)-\text{linear maps } \gamma:V \arrow V \text{ which preserve the form } h          \\
     &  & \text{ and map }
        L[1/p_1,/p_2,\ldots,1/p_l] \text{ to itself }
\end{eqnarray*}
maps into a uniform lattice $\overline{\Gamma_1} < \PGL_3(\Q_{p_1}) \times \PGL_3(\Q_{p_2})
\times \ldots \times \PGL_3(\Q_{p_l})$. We wish to work with subgroups of $\Gamma_1$ defined
by some congruence conditions mod $7$. Let
\begin{equation*}
\pi_7:\Gamma_1(p_1,\ldots,p_l) = G \left( \Z \left[1/p_1,\ldots,1/p_l
        \right] \right) \arrow
        G \left( \frac{\Z \left[1/p_1,\ldots,1/p_l \right]}
        {7\Z \left[ 1/p_1,\ldots,1/p_l\right]} \right) =
         G(\F_7)
\end{equation*}
be the reduction map mod $7$. The Hermitian form $\overline{h}$ becomes degenerate of rank
$1$ over $\F_7$ with a null space $L_1$ spanned by $\zeta - 1$ and $(\zeta -1)^2$.
The action of $\Gamma_1(p_1,\ldots,p_l)$ on $L_1$ gives a homomorphism
$$\eta:\Gamma_1(p_1,\ldots,p_l) \arrow \GL_2(\F_7).$$
Let $\Hat{H} < \GL_2(\F_7)$ be a $2$-Sylow subgroup $H = \hat{H} \cap \SL_2(\F_7)$ and
$\Gamma(p_1,\ldots,p_l) = \phi^{-1}(\Hat{H})$. Mumford proves the following:
\begin{proposition} \label{thm:Mumford} (Mumford \cite{Mumford:Ample})
Using the terminology of the above paragraph the group $\Gamma(2)$ is torsion free and acts
transitively on the vertices of the building of $\PGL_3(\Q_2)$.
\end{proposition}
\noindent {\it{Note: }} Mumford in fact shows that $\eta(\Gamma_1(2)) =  \SL_2(\F_7)$. He
defines $\Gamma(2)$ as a pullback of a $2$-Sylow subgroup of $\SL_2(\F_7)$. Since
$\SL_2(\F_7)$ is normal, this is consistent with our definition. In the course of our proof we
will see that, in fact even in our more general setting $\eta(\Gamma_1) \subset \SL_2(\F_7)$.

\begin{proposition}
\label{prop:HR_lattice} Let $p_1,\ldots,p_l$ be primes all congruent to $1,2$ or $4$ ($\mod
7$), $\Delta_i$ the building of $\PGL_3(\Q_{p_i})$, $\Delta_i^{(1)}$ its one skeleton and $Y =
\Delta_1^{(1)} \times \ldots \times \Delta_1^{(1)}$. Then the group $\Gamma(p_1,\ldots,p_l)$
acts freely and transitively on the vertices of $Y$ and it acts freely on the set of geometric
cubes in $Y$.
\end{proposition}
\begin{proof}
We proceed in several steps:
 \begin{enumerate}
 \item Set $\Gamma = \Gamma(p_1,p_2,\ldots,p_l)$ and $\Gamma_1 =
\Gamma_1(p_1,p_2,\ldots,p_l)$. Assume first that $2$ is one of the given primes, say $p_1 =
2$. We denote by $\proj$ the projection on all but the first factor:
$$\proj:\PGL_3(\Q_{p_1}) \times \ldots \times \PGL_3(\Q_{p_l}) \arrow
      \PGL_3(\Q_{p_2}) \times \ldots \times \PGL_3(\Q_{p_l}).$$
\item {\bf{$\proj(\Gamma_1)$} is transitive on the vertices of $\Delta_2^{(1)} \times \ldots \times
\Delta_l^{(1)}$.} \label{sec:transitive_pr} In fact, by weak approximation $\proj(\Gamma_1)$
intersects with $\PSL_3(\Q_{p_2}) \times \ldots \times \PSL_3(\Q_{p_l})$ in a dense subgroup,
thus it acts transitively on each one of the $3^{l-1}$ different vertex colors. It remains
only to exhibit, for each $p \in \{p_2,p_3,\ldots,p_l\}$, an element $\rho_p \in \Gamma_1$
with
    $$val_{q}(det(g)) \left\{
       \begin{tabular}{ll}
          $\ne 0 \ (\mod 3)$ &  $q=p$ \\
          $= 0 \ (\mod 3)$ & $p \ne q$
       \end{tabular} \right. $$
     By assumption $p = \mathfrak{P} \overline{\mathfrak{P}}$ splits in $K$. Since $K$ is a UFD
we may assume that $\mathfrak{P} = (P) = (s + t\sqrt{-7})$ and $\overline{\mathfrak{P}} =
(\overline{P}) = (s - t\sqrt{-7})$ are principal ideals. We fix the notation in such a way
that in $\Q_p = K_{\mathfrak{P}}$ we have $P = \text{(unit)} \cdot p$ and $\overline{P} =
\text{(unit)}$. Define:
    \begin{equation}
    \label{eqn:rho_p}
    \rho_p =
    \begin{pmatrix}
        1 & 0 &  0      \\
        0 & 1 & 0       \\
        \frac{2}{\overline{P}} \frac{P - \overline{P}}{\sqrt{-7}}  &
        -\frac{\lambda^{2}}{\overline{P}}\frac{P - \overline{P}}{\sqrt{-7}}  &
        \frac{P}{\overline{P}} \\
    \end{pmatrix} =
    \begin{pmatrix}
        1 & 0 & 0         \\
        0 & 1 & 0         \\
        \frac{4t}{\overline{P}} & \frac{-2t\lambda^{2}}{\overline{P}} & \frac{P}{\overline{P}} \\
    \end{pmatrix}
    \end{equation}
\item
 {\bf{Identifying $\Gamma(2)$ in $\Gamma$.}} Let $O = (O_1,O_2,\ldots,O_l) \in \Delta =
\Delta_1 \times \ldots \times \Delta_l$ be a base vertex, we choose the notation so that $O_i
\in \Delta_i$ is the vertex whose stabilizer is $\PGL_3(\Z_{p_i}) =
\PGL_3(\mathcal{O}_{\mathfrak{P}_i})$. We claim that
$$\Gamma_1(2) = \Gamma_1 \cap \left[\GL_3(\Q_{p_1}) \times \GL_3(\Z_{p_2}) \times \ldots \times
\GL_3(\Z_{p_l})\right].$$ In fact, any element of the right hand side automatically has all
its matrix elements in $\Z[\lambda][1/2,1/\overline{P_2},\ldots,1/\overline{P_l}]$. Since
$^t\overline{g}=H^{-1}g^{-1}H$, the matrix coefficients of $^t\overline{g}$ are also in the
same ring so that in fact all matrix coefficients are in $\Z[\lambda][1/2]$. The other
inclusion is clear.

We can therefore identify $\Gamma_1(2) = (\Gamma_1)_{(O_2,O_3,\ldots,O_l)}$ where the action
of $\Gamma_1$ on the right hand side is via the projection $\proj$. The action
$(\Gamma_1)_{(O_2,O_3,\ldots,O_l)} \leadsto \Delta_1 \times \{(O_2,O_3,\ldots,O_l)\}$ is
isomorphic to the action $\Gamma_1(2) \leadsto \Delta_1$. Furthermore reduction $\mod 7$ is
compatible with the embedding of $\Gamma_1(2) \hookrightarrow \Gamma_1$: restricting $\eta$ to
$(\Gamma_1)_{(O_2,O_3,\ldots,O_l)}$ gives the reduction map mod $7$ on $\Gamma_1(2)$ . From
here on we identify $\Gamma_1(2)$ with a subgroup of $\Gamma_1$ and $\Gamma(2)$ with a
subgroup of $\Gamma$.
\item {\bf{$\eta(\Gamma_1) \subset \SL_2(\F_7)$.}}
We will in fact prove $\Gamma_1 = \langle \Gamma_1(2),\rho_{p_1},\rho_{p_2},\ldots,\rho_{p_l}
\rangle$ where $\rho_{p_i}$ are the group elements given in equation \ref{eqn:rho_p}. This
will complete the proof because it is already proved by Mumford that $\eta(\Gamma_1(2))
\subset \SL_2(\F_7)$, and for the elements $\rho_p$ it can easily be verified\footnote{ In
Mumfords notation \cite{Mumford:Ample} the matrix is transposed to this one} that in the basis
$1,\zeta - 1 ,(\zeta-1)^{2}$:
\begin{equation}
\eta(\rho_p) =
\begin{pmatrix}
    1 & 0 & 0 \\
    0 & 1 & 0 \\
    0 & \frac{5}{\overline{P}} \frac{P - \overline{P}}{\sqrt{-7}} & 1 \\
\end{pmatrix}.
\end{equation}

Set $\Gamma_1' = \langle  \Gamma_1(2),\rho_{p_1},\rho_{p_2},\ldots,\rho_{p_l} \rangle$, both
$\Gamma_1$ and $\Gamma_1'$ have the same vertex stabilizers acting on the product $\Delta_2
\times \Delta_3 \times \ldots \times \Delta_l$, namely $(\Gamma_1')_{(O_2,O_3,\ldots,O_l)} =
(\Gamma_1)_{(O_2,O_3,\ldots,O_l)} = \Gamma_1(2)$. In order to show that they are isomorphic we
only have to show that $\Gamma_1'$ acts transitively on the vertices of $\Delta_2 \times
\ldots \times \Delta_l$. The vertex $(O_2,O_3,\ldots,O_l) \in \Delta_2 \times \ldots \times
\Delta_l$ has $2^{l-1}$ different types of vertices among its nearest neighbors, and we can
reach at least one representative of each of these types using the elements
$\rho_2,\rho_3,\ldots,\rho_l$ and their adjoints. By weak approximation theorem $\Gamma_1(2)$
acts transitively on the nearest neighbors of the same type, so we can reach all the nearest
neighbors of $(O_2,,O_3,\ldots,O_l)$. This suffices to prove transitivity.

\item
{\bf{$\Gamma_1$ acts transitively on the vertices of $Y$.}} \label{sec:transitive_1} Let $V =
(V_1,V_2,\ldots,V_l)$ be a vertex of $Y$. By \ref{sec:transitive_pr} there is an element
$\gamma \in \Gamma_1$ such that $\gamma V_i = O_i \ \forall \ 2 \le i \le l$. By Mumfords'
theorem (proposition \ref{thm:Mumford} above) and by the identification $\Gamma_1(2) =
(\Gamma_1)_{(O_2,\ldots,O_l)}$ we can find an element $\delta \in
(\Gamma_1)_{(O_2,\ldots,O_l)}$ such that $\delta \gamma V = O$.
\item
{\bf{$\Gamma$ acts freely.}} To prove that $\Gamma$ acts freely on the geometric cubes it is
enough to show that $\Gamma$ acts freely on the vertices, because even the ambient group
$\PGL_3(\Q_p)$ does not invert edges in $\Delta_p^{(1)}$ (note that this is no longer true for
$2$-simplexes in $\Delta_p$).

By Mumfords' theorem (proposition \ref{thm:Mumford}) the $\Gamma$ stabilizer of $O$ is
trivial:
\begin{eqnarray*}
\Gamma_{O} & = & \left(\Gamma_{(O_2,O_3,\ldots,O_l)}\right)_{O_1}       \\
           & = & \left( \Gamma(2) \right)_{O_1}                         \\
           & = & \langle e \rangle
\end{eqnarray*}
Since Mumfords' proposition is true for any choice of a $2$-Sylow subgroup we know even more:
$\eta((\Gamma_1)_O)$ intersects trivially any $2$-Sylow subgroup of $\hat{G}$.
This together with the transitivity of the action of $\Gamma_1$
(proposition\ref{sec:transitive_1})
proves that every vertex stabilizer is trivial. In fact, given another vertex $V$,
there is a $\gamma \in \Gamma_1$ such that $\gamma O = V$,
and by our observation above $\eta((\Gamma_1)_V) =\eta(g)\eta((\Gamma_1)_O) \eta(g^{-1})$
intersects trivially with $\hat{H}$.
\item
 {\bf{$\Gamma$ acts transitively.}} The number of orbits for the action of $\Gamma$ on
the vertices of $Y$ is exactly $\tfrac{[\Gamma_1:\Gamma]}{|(\Gamma_1)_{O}|}$. Similarly by
Mumfords' proposition we have $\tfrac{[\Gamma_1(2):\Gamma(2)]}{|(\Gamma_1(2))_{O_1}|}= 1$.
Dividing and using the equality of the denominators we obtain:
\begin{eqnarray*}
  \text{ \# of orbits }
       & = & \frac{[\Gamma_1:\Gamma]}{[\Gamma_1(2):\Gamma(2)]} \\
       & = & \tfrac{[\eta(\Gamma_1):\eta(\Gamma)]}{[\eta(\Gamma_1(2)):\eta(\Gamma(2))]} \\
       & = & \tfrac{[\SL(\F_7):H]}{[\SL_2(\F_7):H]} = 1
\end{eqnarray*}
Where in the third equality we have used the fact (see \cite{Mumford:Ample}) that
$\eta(\Gamma_1(2)) = \SL_2(\F_7)$ and consequently $\eta(\Gamma(2)) = H$ $2$-Sylow subgroup.
 \item
 {\bf{The general case.}} If $2 \not \in \{p_1,\ldots,p_l\}$ then we can add it
artificially by setting $p_0 = 2$. Now by identifying the two actions:
$\Gamma(p_1,p_2,\ldots,p_l) \action \Delta_1 \times \Delta_2 \times \ldots \times \Delta_l$
and $(\Gamma(p_0,p_1,\ldots,p_l))_{O_0} \leadsto \{O_{0}\} \times \Delta_1 \times \Delta_2
\times \ldots \times \Delta_l$ we deduce that the latter action satisfies the desired
properties. This completes the proof of proposition \ref{prop:HR_lattice}
\end{enumerate}
\end{proof}
\subsubsection{A Kazhdan automaton}
The example below is a particular example of a lattice acting on a product of trees, coming
from an arithmetic lattice in higher rank $p$-adic groups. These lattices play an important
role in the in the classification of lattices in products of trees that admit infinite,
finite-dimensional representations in characteristic $0$ \cite{BMZ:Representations}. A similar
construction appears also in \cite{BM:Lattices}.

\begin{proposition}
\label{prop:T-automata} Let $\Gamma < \PGL_3(\Q_p) \times \PGL_3(\Q_q)$ be an irreducible
lattice that acts vertex transitively and freely on the (geometric cells of the) product of
the one skeletons of the buildings $Y = \Delta_p^{(1)} \times \Delta_q^{(1)}$. Then $\Sigma =
\Gamma \bs Y$ is a VH-T-square complex with one vertex and $G(A^{T}(\Sigma)) = \Gamma_{O_p}$
is a lattice in $\PGL_3(\Q_q)$. In particular it has property $T$.
\end{proposition}
\noindent
{\it{Note: }} The situation is symmetric, so by changing $p$ and $q$ we can show that
$G(A^{T}(\Sigma)^{*}) = \Gamma_{O_q}$ is a lattice in $\PGL_3(\Q_p)$ and has property
$T$ as well.
\begin{proof}
It is clear that $\Sigma$ is a VH-T-complex and therefore its universal cover is a product
of trees $\tilde{\Sigma} = \tilde{Y} = T_h \times T_v$, where $T_h$ is a $2(p^2+p+1)$-regular
tree and $T_v$ is a $2(q^2+q+1)$-regular tree. In fact we obtain the following tower of
regular coverings:
\begin{equation*}
 (T_h \times T_v,\overline{O})  \arrow
 (Y,\tilde{O} ) \arrow
 (\Sigma,O)
\end{equation*}
Here $\overline{O} = (\overline{O}_1,\overline{O}_2)$ and $\tilde{O} = (\tilde{O}_1,\tilde{O}_2)$
are base vertices.
The covering groups fit into a short exact sequence:
\begin{equation*}
   1 \arrow \pi_1(Y,\tilde{O})
     \arrow \Theta \defeq \pi_1(\Sigma,O)
     \stackrel{p}{\arrow} \Gamma \arrow 1,
\end{equation*}
where we have defined $\Theta \defeq \pi_1(\Sigma,O)$.

By proposition \ref{prop:action_description_T} we are interested in the group $\left( \proj_1
(\Theta) \right)_{\overline{O_1}}$, where $\proj_1: \Aut(T_h) \times \Aut(T_v) \arrow
\Aut(T_h)$ is just the map that forgets the action on the $T_v$ coordinate. First let us note
the following
\begin{lemma}
\begin{equation*}
   \ker(\proj_1) \cap \Theta = \pi_1(\Delta_1^{(1)},\tilde{O_2})
\end{equation*}
\end{lemma}
\begin{proof}
Indeed assume that $\theta \in \ker(\proj_1) \cap \Theta$, then $p(\theta) = (e,\rho) \in
\ker(\proj_1) \cap \Gamma$, where by abuse of notation we use the same notation for $\proj_1$
downstairs. By weak approximation one can find elements $\gamma=(\gamma_1,\gamma_2) \in
\Gamma$ with $\gamma_2$ arbitrarily close, but not equal to $\rho$ and therefore:
\begin{equation*}
[\gamma,p(\theta)] = [(1,[\gamma_2,\rho])] \stackrel{\gamma_2 \arrow \rho}{\longrightarrow} e
\end{equation*}
Since $\Gamma$ is discrete this implies that $p(\theta) = e$, and $\rho \in
\pi_1(\Delta_1^{(1)},\tilde{O_2})$. Since the opposite inclusion is obvious this concludes the
proof of the lemma.
\end{proof}

Since $\pi_1(\Delta_p^{(1)})$ acts freely on $T_h$ we have:
\begin{eqnarray*}
  \left(\proj_1 (\Theta) \right)_{\overline{O_1}}
  & = & \left( p (\proj_1(\Theta))\right)_{\tilde{O_1}} \\
  & = & p \left( \frac{\Theta}{\pi_1(\Delta_p^{(1)},\tilde{O_1})} \right)_{\tilde{O_1}} \\
  & = & \left( \frac{\Theta}  {\pi_1(\Delta_p^{(1)},\tilde{O_1})
                               \pi_1(\Delta_q^{(1)},\tilde{O_2})} \right)_{\tilde{O_1}}  \\
  & = & \Gamma_{\tilde{O_1}}
\end{eqnarray*}
completing the proof of proposition \ref{prop:T-automata}.
\end{proof}
\noindent
\bibliography{yair}
\end{document}